\documentclass[notitlepage, 10pt]{article}
\title{Aztec Diamonds and Baxter Permutations}
\author{{\sc Hal Canary} \\ {\small University of Wisconsin--Madison}}
\date{March 2, 2004}
\usepackage{amssymb}
\usepackage{hyperref}
\usepackage{graphicx}
\hypersetup{pdfstartview=FitB,pdfauthor={Hal Canary},
pdftitle={Aztec Diamonds and Baxter Permutations}}

\newtheorem{theorem}{Theorem}
\newtheorem{lemma}{Lemma}

\begin{document}

\maketitle

\begin{abstract}
We present a proof of a conjecture about the relationship between Baxter
permutations and pairs of alternating sign matrices that are produced
from domino tilings of Aztec diamonds.  It is shown that if and only
if a tiling corresponds to a pair of ASMs that are both permutation
matrices, the larger permutation matrix corresponds to a Baxter
permutation.

There has been a thriving literature on both pattern-avoiding
permutations of various kinds \cite{baxter} \cite{dulucq} and tilings
of regions using dominos or rhombuses as tiles \cite{EKLP} \cite{kuo}.
However, there have not as of yet been many links between these two
areas of enumerative combinatorics.  This paper gives one such link.
\end{abstract}

\section{Introduction}

\begin{figure*}[h]
\begin{center}
\includegraphics{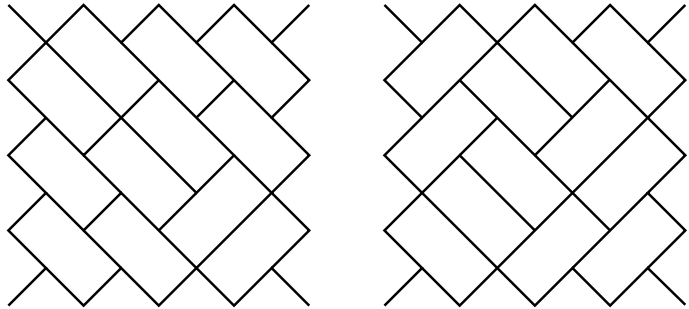}
\caption{ \small \label{toad} Two examples of tilings of Aztec
diamonds.}
\end{center}
\end{figure*}

Figure \ref{toad} shows two of the 64 different ways of tiling an
order three Aztec diamond with dominoes.  Aztec Diamonds are defined
in ``Alter\-nating-Sign Matrices and Domino Tilings,'' by Elkies,
Kuperberg, Larsen, and Propp \cite{EKLP}.  Each tiling will be
referred to as a TOAD, a Tiling Of an Aztec Diamond.  In all of the
illustrations, we have rotated the Aztec Diamonds and drawn spurs on
each corner to simplify later calculations.

An Alternating Sign Matrix (ASM) is a square matrix whose entries are
$0$, $1$, or $-1$, such that the entries of each row or column must
add up to $1$, and all nonzero entries must alternate in sign.
Examples:
\begin{displaymath}
\left( \begin{array}{ccc}
0 & 1 & 0 \\
0 & 0 & 1 \\
1 & 0 & 0 
\end{array} \right) \qquad
\left( \begin{array}{cccc}
 0  &  1  &  0  &  0  \\
 0  &  0  &  0  &  1  \\
 1  & -1  &  1  &  0  \\
 0  &  1  &  0  &  0
\end{array} \right) 
\end{displaymath}

Elkies, Kuperberg, Larsen, and Propp showed that there are
$2^{n(n+1)/2}$ tilings of an order $n$ Aztec diamond with dominoes
\cite{EKLP}.  As part of their proof, they used a relationship between
TOADs and pairs of comatable ASMs.
\emph{Compatibility} of ASMs was defined by Robbins and Rumsey in
terms of the ASMs themselves, without reference to domino tilings
\cite{robbins}.

\begin{figure*}[h]
\begin{center}
\includegraphics{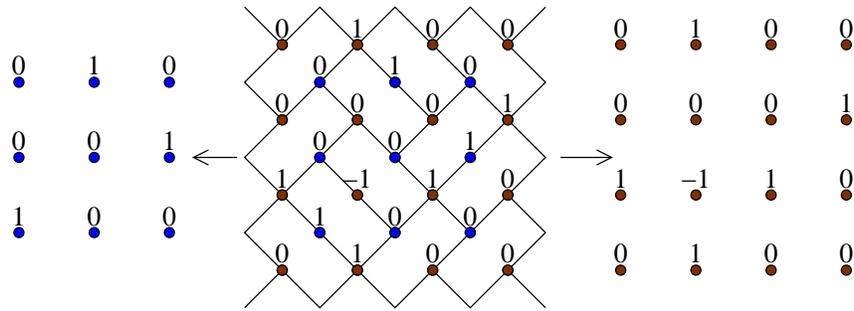}
\caption{ \small \label{asmtoad} How to produce two ASMs from a TOAD.}
\end{center}
\end{figure*}

The interior vertices of an order $n$ TOAD
are arranged as a $n$-by-$n$ square matrix inside an
$(n+1)$-by-$(n+1)$ matrix.  Each interior vertex is connected to two
or more of its four nearest neighbors.
For the smaller square assign a
$0$ to each vertex incident to exactly three edges, a $1$ to each
vertex incident to exactly two edges, and a $-1$ to each vertex
incident to exactly four edges.  For the larger matrix, do the same thing,
but reverse the roles of $1$'s and $-1$'s. The two above examples of
ASMs come from the Aztec Diamond in Figure \ref{asmtoad}.

A pair of ASMs of order $n$ and order $n+1$ is said to be
compatible if there is a TOAD that is related to the pair by
the above algorithm.   A pair of compatible ASMs uniquely determines a
TOAD.  
It has been shown that an order $n+1$ ASM with $k$ $-1$'s is
compatible with $2^k$ order $n$ ASMs.  Consequently an ASM with no
$-1$'s is only compatible with one smaller ASM.  An ASM with no $-1$'s
is a permutation matrix. 

The study of Baxter permutations began outside of the field of
combinatorics \cite{baxter}, and even though they show up many places,
it was surprising to find a connection between domino tilings and
Baxter permutations.  We will give two definitions of Baxter
permutations.  The latter follows from the former and is described in
terms of permutation matrices.

A permutation $\pi$ of $n$ objects is a Baxter permutation if for each 
$i \in \{1,2,3,\ldots ,n-1\}$,  there is a $k_i$ such that
$\pi(k_i)$ is between $\pi(i)$ and $\pi(i+1)$ (allowing $k_i$
to be $i$ but not $i+1$) and for every
$\pi(m)$ between $\pi(i)$ and $\pi(k_i)$ 
(including $\pi(k_i)$), $m\leq i$, 
and for every $\pi(m)$ between $\pi(k_i)$ and $\pi(i+1)$ 
(not including $\pi(k_i)$), $m > i+1$ \cite{dulucq}.

The easiest way to understand this definition is by looking at
at a permutation matrix. Let $B$ be a permutation matrix.
By definition, each row and each column of $B$ contains exactly
one entry that is a 1:  for each $i$, $B_{i,\pi(i)} = 1$ 
and the rest are zeros.

Take any two adjacent rows in $B$, the $i$th and $(i+1)$th.
Then $B_{i,\pi(i)} = B_{i+1,\pi(i+1)} =1$.  Then, we look at the 
columns between the $\pi(i)$th column and the $\pi(i+1)$th 
column.

If the permutation is Baxter then there will be a vertical dividing
line between two of these columns, so that every 1 on the
$\pi(i+1)$ side of the line is in a row below the $i+1$th row, and
every 1 on the $\pi(i)$ side of the line is in a row above the
$i$th row.  For example, let us test permutation 45123  between the
second and third rows:  
\begin{displaymath}
B =
\left( \begin{array}{ccc|cc}
 0  &  0  &  0  &  1  &  0  \\
 0  &  0  &  0  &  0  &  1  \\
\hline
 1  &  0  &  0  &  0  &  0  \\
 0  &  1  &  0  &  0  &  0  \\
 0  &  0  &  1  &  0  &  0 
\end{array} \right)
\end{displaymath}
Because there is a place to draw a vertical line between columns 1 and
5, this matrix passes the test for $i=2$.  Since it passes the test
for each of the other rows as well, it is Baxter. 

Now we are ready to present our theorem, which we will prove in
Section \ref{proof}:

\begin{theorem}\label{bpth}
An order $(n+1)$ ASM without $-1$'s is compatible with an
order $n$ ASM without $-1$'s if and only if it is a Baxter
permutation matrix.  
\end{theorem}

As a consequence of this theorem, we know the number of order $n$
TOADs with the property that vertices corresponding to the smaller ASM
never have four edges incident to them, and vertices corresponding to
the larger ASM never have two edges incident to them.  That number is
the same as the number of Baxter permutations of order $n+1$, which
is given by Chung, Graham, Hoggatt, and Kleiman \cite{count}.  The
number of Baxter permutations on order $n$ is 
 \begin{displaymath}
 \sum_{r=0}^{n} 
 \frac{{n+2 \choose r}{n+2 \choose r+1}{n+2 \choose r+2}}
 {{n+2 \choose 1}{n+2 \choose 2}}.
 \end{displaymath}

\section{Construction of the smaller ASM}\label{sasm}

To prove Theorem \ref{bpth}, we will present an algorithm that will
produce the smaller ASM that is compatible with a given permutation
matrix.  The following lemma is a formal description of the algorithm.
We will use the relationship between ASMs and Aztec diamonds to prove
the lemma by the end of section \ref{sasm}.

\begin{lemma}\label{alg}
If $B$ is an order $n+1$ permutation matrix, and $A$ is the order $n$
ASM that is compatible with $B$, then $A_{i,j} = 0$ if and only if 
\center
$(\forall~ k \leq i)  ~ B_{k,j} = B_{k,j+1} = 0$ 
or $(\forall~ k \geq i+1)~ B_{k,j} = B_{k,j+1} = 0$ 
or $(\forall~ k \leq j)  ~ B_{i,k} = B_{i+1,k} = 0$ 
or $(\forall~ k \geq j+1)~ B_{i,k} = B_{i+1,k} = 0$.
\flushleft
The nonzero elements of $A$ alternate between  $1$ and $-1$.
\end{lemma}

\subsection{Using the Algorithm}

We will use an example to show how to apply the algorithm.  Suppose we
are given the permutation 31425.  First, let us define the matrix $B$
which corresponds to our permutation.  

\begin{displaymath}
B = 
\left( \begin{array}{ccccc}
 0  &  0  &  1  &  0  &  0  \\
 1  &  0  &  0  &  0  &  0  \\
 0  &  0  &  0  &  1  &  0  \\
 0  &  1  &  0  &  0  &  0  \\
 0  &  0  &  0  &  0  &  1 
\end{array} \right)
\end{displaymath}

Now draw the smaller matrix $A$ interspersed within $B$.
\begin{displaymath}
\left( \begin{array}{ccccccccc}
 0  &        &  0  &        &  1  &        &  0  &        &  0  \\
    & A_{11} &     & A_{12} &     & A_{13} &     & A_{14} &     \\
 1  &        &  0  &        &  0  &        &  0  &        &  0  \\
    & A_{21} &     & A_{22} &     & A_{23} &     & A_{24} &     \\
 0  &        &  0  &        &  0  &        &  1  &        &  0  \\
    & A_{31} &     & A_{32} &     & A_{33} &     & A_{34} &     \\
 0  &        &  1  &        &  0  &        &  0  &        &  0  \\
    & A_{41} &     & A_{42} &     & A_{43} &     & A_{44} &     \\
 0  &        &  0  &        &  0  &        &  0  &        &  1 
\end{array} \right)
\end{displaymath}

From the description of the algorithm in Lemma \ref{alg}, 
$A_{11}=0$ because $B_{11}=B_{12}=0$ and 
$A_{14}=0$ because $B_{14}=B_{15}=0$ and 
$A_{24}=0$ because $B_{14}=B_{15}=B_{24}=B_{25}=0$. 
Repeated applications of this method will give :
$ A_{11} = A_{14} = A_{13} = A_{24} = A_{31} = A_{33} $
$= A_{34} = A_{41} = A_{42} = A_{43} = 0.$ 
Graphically:

\begin{displaymath}
\left( \begin{array}{ccccccccc}
 0  &        &  0  &        &  1  &        &  0  &        &  0  \\
    &      0 &     & A_{12} &     &      0 &     &      0 &     \\
 1  &        &  0  &        &  0  &        &  0  &        &  0  \\
    & A_{21} &     & A_{22} &     & A_{23} &     &      0 &     \\
 0  &        &  0  &        &  0  &        &  1  &        &  0  \\
    &      0 &     & A_{32} &     &      0 &     &      0 &     \\
 0  &        &  1  &        &  0  &        &  0  &        &  0  \\
    &      0 &     &      0 &     &      0 &     & A_{44} &     \\
 0  &        &  0  &        &  0  &        &  0  &        &  1 
\end{array} \right)
\end{displaymath}

Lemma \ref{alg} tells us that the other elements---$A_{12}$,
$A_{21}$, $A_{22}$, $A_{23}$, $A_{32}$, and $A_{44}$---are nonzero and
alternate between 1 and $-1$.  Thus,
\[A = \left( \begin{array}{cccc}
0	& 1	& 0	& 0\\
1	& -1	& 1	& 0\\
0	& 1	& 0	& 0\\
0	& 0	& 0	& 1\\
\end{array} \right)\]

The proof of Lemma \ref{alg} comes in two parts.  First we show that
entries in $A$ that have a certain property are equal to 0.  Second, we
show that entries that lack this property are nonzero.  The values of
the nonzero entries are given by the definition of an ASM.

\subsection{Zeros}

First let us show that  $(\forall~ k \leq j)~(B_{i,k} = B_{i+1,k} =
0 \Rightarrow A_{i,j}=0)$ by induction on $j$.  

\begin{figure*}[h]
	  \begin{center}
	    \includegraphics{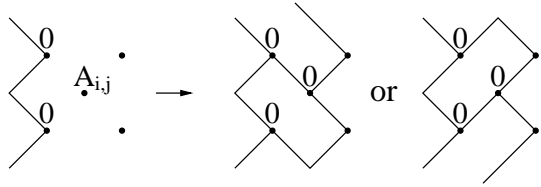}
	    \caption{ \small \label{varprop} Determining $A_{i,1}$ }
	  \end{center}
	\end{figure*}

Let $j=1$.  Figure \ref{varprop} shows a situation where $B_{i,1} =
B_{i+1,1}=0$.  We need to show that $A_{i,1}$ must be zero.  To do
this, we look at the positioning of the dominoes in the TOAD that
corresponds to $A$ and $B$.  The square located between the vertices
of $B_{i,1}$ and $B_{i+1,1}$ must be part of a domino.  Since it is on
the left edge of the Aztec diamond, the domino must be pointing
northeast or southeast.  
Figure \ref{varprop} shows both cases. The figure also shows the
$B_{i,1}$ and $B_{i+1,1}$ vertices with three edges incident to them.  
Either configuration leads to $A_{i,1}=0$.

For the induction step choose any $j \ge 1$ such that all of the
$B_{i,k}$'s and $B_{i+1,k}$'s to the left of $A_{i,j}$ are zero.
We assume not only that all of the $A_{k,j}$'s to the left of
$A_{i,j}$ are zero, but also that the zig-zag shape from the far left
edge of the Aztec diamond has propagated inward towards the $j$th
column. 
	\begin{figure*}[h]
	  \begin{center}
	    \includegraphics{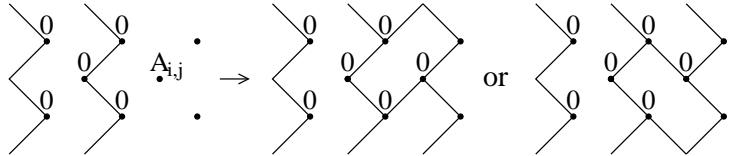}
	    \caption{ \small \label{induction}The induction step.  }
	  \end{center}
	\end{figure*}
Since $B_{i,j}=B_{i+1,j}=0$, both of these vertices meets exactly
three edges. 
As Figure \ref{induction} shows, once we know the shape of the tiling
to the left of $A_{i,j}$, and that $B_{i,j}=B_{i+1,j}=0$, we can
determine the value of $A_{i,j}$ just as we determined the value of
$A_{i,1}$. 

\begin{figure*}[h]
	  \begin{center}
	    \includegraphics{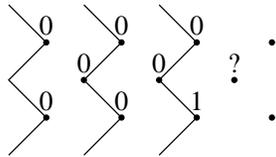}
	    \caption{ \small \label{prop} 
	      Zeros in $A$ propagate until they hit a 1 in $B$.}
	  \end{center}
	\end{figure*}

Figure \ref{prop} shows that this induction can continue until we
encounter a 1.

In other words, 0's propagate from the edge of $A$ until they hit a
nonzero entry in the surrounding rows and columns of $B$.
This proof can also be applied to the top, bottom, and right edges of
the ASMs.  Then we have proved half of the lemma.  

\subsection{ Non-zero matrix entries }

To finish proving Lemma \ref{alg}, we will need to look at how the
tiling patterns must propagate inside the matrix $A$.  
We will concentrate on those elements of $A$ whose value is still
undetermined.  To represent those entries, we will use the symbol
$\Box$.

We will now look at one individual row of $A$, and look at the two
rows of $B$ that are above and below it.
Here is an example: the second and third rows from the permutation 34215.
\begin{displaymath}
\begin{array}{ccccccccc}
 0  &     &  0  &     &  0  &     &  1  &     &  0  \\
    &  0  &     & \Box&     &  0  &     &  0  &     \\
 0  &     &  1  &     &  0  &     &  0  &     &  0  \\
\end{array}
\end{displaymath}
Notice first that all of the unknown entries on any row in $A$ will
be between the 1's in the two closest rows of $B$.
As Figure \ref{direction} shows, a 0 will show up between the 1's only if
there is a vertical path of zeros from
either the top edge or the bottom edge. In this example, the zero is
there because of a path from the bottom edge.  
We can say that a $\Box$ arises when the columns of $B$ to the
immediate
left and right of the box have their 1's in opposite directions, and
that if these two columns of $B$ have their 1's in the same direction,
there is a 0.  

	\begin{figure}[h]
\begin{displaymath}
\begin{array}{ccccccccc}
    &     &     &     &  1  &     &    &     &     \\
    &     &     &     & \vdots &  &  &  &     \\
 0  &     &  0  &     &  0  &     &  1  &     &  0  \\
    &  0  &     &\Box &     &  0  &     &  0  &     \\
 0  &     &  1  &     &  0  &     &  0  &     &  0  
\end{array}
\end{displaymath}
	\caption{\small \label{direction} From the position of the
	$\Box$s in a row of $A$, we can reconstruct where 1's are
	located in $B$.  }
	\end{figure}

Our next example is the (non-Baxter) permutation 3142.  We choose this
because the second row of $A$ will have three $\Box$s.  
\begin{displaymath}
\begin{array}{ccccccc}
  0  &     &  0  &     &  1  &     &  0       \\
     &  0  &     & \Box&     &  0  &          \\
  1  &     &  0  &     &  0  &     &  0       \\
     &\Box &     & \Box&     & \Box&          \\
  0  &     &  0  &     &  0  &     &  1       \\
     &  0  &     & \Box&     &  0  &          \\
  0  &     &  1  &     &  0  &     &  0     
\end{array}
\end{displaymath}
Looking only at the center row, we can reconstruct where the 1's are
in $B$:
\begin{displaymath}
\begin{array}{ccccccc}
     &     &     &     &  1  &     &        \\
     &     &     &     &\vdots &   &        \\
  1  &     &  0  &     &  0  &     &  0     \\
     &\Box &     & \Box&     & \Box&        \\
  0  &     &  0  &     &  0  &     &  1     \\
     &     & \vdots &  &     &     &	   \\
     &     &  1  &     &     &     &
\end{array}
\end{displaymath}
It can now be seen that there is an odd number of boxes in each row,
because each box always corresponds to a switching in the direction of
the 1's, and the 1 on the left of the leftmost box must be 
in the direction opposite of the 1 in the column to the right of the
rightmost box. 

Consequently, $\Box$s can be in either of two configurations.  
The important thing is the direction of the 1 in
the two adjacent rows and two adjacent columns of $B$.

There are only two types of configurations, up to various symmetries,
as Figure \ref{windmill} illustrates.
Call the configuration represented by the one on the right
\emph{windmilled}, and the configuration represented by the one on the
left \emph{non-windmilled}.

	\begin{figure}[h]
\[
\begin{array}{ccccccc}
	&	& 0	&	& 1	&	&	\\
	&	&\vdots	&	&\vdots	&	&	\\
0	&\cdots	& 0	&	&  	&\cdots & 1	\\
 	&	&	&\Box	&	&	&	\\
1	&\cdots	&  	&	& 0	&\cdots & 0	\\
	&	&\vdots	&	&\vdots	&	&	\\
	&	& 1	&	& 0	&	& 	\\
\end{array}
\quad \textrm{or} \quad
\begin{array}{ccccccc}
	&	& 1	&	& 0	&	&	\\
	&	&\vdots	&	&\vdots	&	&	\\
0	&\cdots	& 0	&	& 0	&\cdots & 1	\\
 	&	&	&\Box	&	&	&	\\
1	&\cdots	& 0	&	& 0	&\cdots & 0	\\
	&	&\vdots	&	&\vdots	&	&	\\
	&	& 0	&	& 1	&	& 	\\
\end{array}
\]
	\caption{\small \label{windmill} The non-windmilled (left) and
	windmilled (right) $\Box$ configurations.}
	\end{figure}

It can be seen that the $\Box$s along a row alternate between
non-windmilled and  windmilled, starting with a non-windmilled $\Box$.
It is claimed that the left configuration  will lead to $\Box = 1$
and the  windmilled always lead to $\Box = -1$. To prove this, we will
have to look at the tiling patterns between the 1's in two consecutive
rows of $B$.

\subsection{Tilings and $\Box$s}

In this section, we will assume for a moment that we know the values
of each entry in the $i$th row of $A$ and the $i$th and $i+1$th rows
of $B$.  
The reader might find it useful to generate example rows and try to
draw the corresponding tilings. 
For example, these rows: 
\[\begin{array}{ccccccccccccc}
0 &   & 0 &   & 0 &   & 0 &   & 0 &   & 1 &   & 0 \\
  & 0 &   & 1 &   & 0 &   & -1&   & 1 &   & 0     \\
0 &   & 1 &   & 0 &   & 0 &   & 0 &   & 0 &   & 0
\end{array}
\]
will produce this partial tiling:
	  \begin{center}
	    \includegraphics{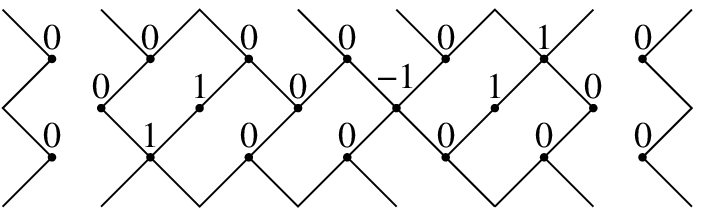}
	  \end{center}
Notice that in the interesting area---between the 1's in $B$---the
tiling is uniquely determined by these three matrix rows.  We will
prove this, and use the machinery of the proof to show that $\Box$s
are always non-zero.
	\begin{figure}[h]
	  \begin{center}
	    \includegraphics{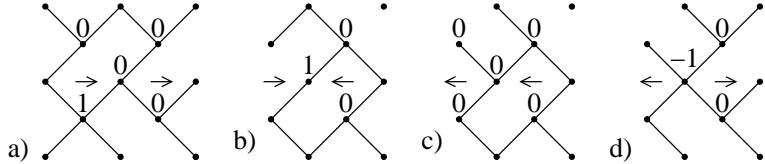}
	    \caption{ \small \label{fill} 
	      Here are examples of how the titling patterns propagate
	      from left to right. 
	      The arrows serve to highlight the
	      change that $1$'s and $-1$'s make to the pattern.
	      \textbf{(a)} and \textbf{(c)} show configurations that
	      go unchanged when they encounter a 0 in $A$.
	      \textbf{(b)} and \textbf{(d)} show how encountering a 1
	      or a $-1$ in $A$ changes the configuration.
	      }
	  \end{center}
	\end{figure}

Figure \ref{fill} shows how the tiling patterns that propagate in
from the left side of the TOAD change first when they hit a 1 in 
$B$ and later when they hit 1's and $-1$'s in $A$.  

\begin{figure}[h]
	  \begin{center}
	    \includegraphics{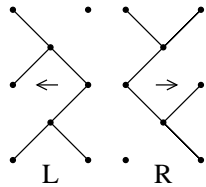}
	    \caption{ \small \label{RL} Once we are to the right of
	      a 1 in $B$, the tiling can have two patterns, $L$ and
	      $R$.  $L$ corresponds to having the closest non-zero $A$
	      entry on the left being 1, and $R$ corresponds to being on
	      the right of a 1 in $A$.  If $R$ and $L$ fit together
	      one way, there is a $-1$ in between.  If they are fit
	      together the other way, there is a $1$ in between.  }
	  \end{center}
	\end{figure}

Notice that
in between the two 1's in $B$, the tiling can take on one of two
patterns:  one will happen if the next nonzero element to the right in
$A$ is $1$, and the other will happen if if the next nonzero element
to the right  is $-1$.  Figure \ref{RL} shows the two patterns. 

To help us keep this straight, the components in
Figures \ref{fill} and \ref{RL} have an arrow that points to the right
if the next nonzero element on the right in $A$ is 1, and to the left
if the next nonzero element on the left in $A$ is a 1.

For example, look at this row from an ASM, with the arrows filled in. 
\begin{displaymath}
\begin{array}{ccccccccccccccc}
\rightarrow & 0 & \rightarrow & 1 & \leftarrow & -1 & \rightarrow & 0 &
\rightarrow & 0 & \rightarrow & 1 & \leftarrow & 0 & \leftarrow 
\end{array}
\end{displaymath}

Suppose that this row had been determined by a larger permutation
matrix it was compatible with.  Here is what might look like:

\begin{displaymath}
\begin{array}{ccccccccccccccc}
1 &  & 0 &  & 0 &  & 0 &  & 0 &  & 0 &  & 0 &  & 0  \\
\rightarrow & 0 & \rightarrow & 1 & \leftarrow & -1 & \rightarrow & 0 &
\rightarrow & 0 & \rightarrow & 1 & \leftarrow & 0 & \leftarrow \\
0 &  & 0 &  & 0 &  & 0 &  & 0 &  & 0 &  & 1 &  & 0
\end{array}
\end{displaymath}
	\begin{center}
	  \includegraphics{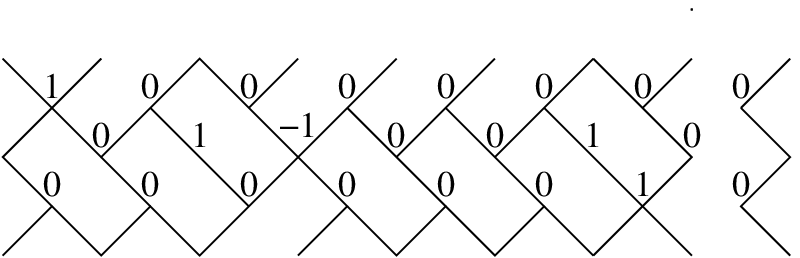}\\
	\end{center}
	\begin{center}
	  \includegraphics{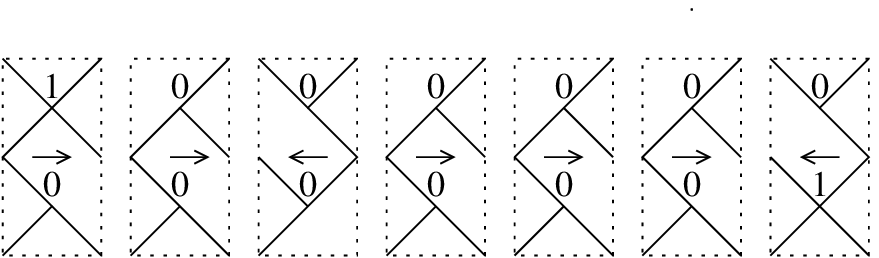}
	\end{center}
The last picture shows the tiling broken into components.  
Notice that the
components between the 1's in  $B$ are all the same, up to a
rotation.  Each rotation corresponds to a nonzero element in the
ASM $A$ that fits between the two components.

This constitutes an algorithm to determine the shape of the local
tiling from the two ASMs.

Now we have assembled all of the tools necessary to finish proving
Lemma \ref{alg}.  We know that the $\Box$s alternate between two
configurations: windmilled and non-windmilled; if we were to show that
the windmilled boxes can not contain 1's or 0's, then the
sign-alternation property of $A$ would force the non-windmilled boxes
to contain 1's, and the lemma would be proved.

\subsection{Windmilled $\Box$s}

Here is what a windmilled configuration looks like, up to reflection.

\begin{displaymath}
\begin{array}{ccccccc}
	&	& 1	&	& 0	&	&	\\
	&	&\vdots	&	&\vdots	&	&	\\
0	&\cdots	& 0	&	& 0	&\cdots & 1	\\
 	&	&	&\Box	&	&	&	\\
1	&\cdots	& 0	&	& 0	&\cdots & 0	\\
	&	&\vdots	&	&\vdots	&	&	\\
	&	& 0	&	& 1	&	&
\end{array}
\end{displaymath}

Let us assume that $\Box =1$, and find a contradiction.  The
contradiction arises when we try to draw the tiling components that
correspond to the arrows that point towards  the 1.
\begin{displaymath}
\begin{array}{ccccccc}
	&	& 1	&	& 0	&	&	\\
	&	&\vdots	&	&\vdots	&	&	\\
0	&\cdots	& 0	& \downarrow& 0	&\cdots & 1	\\
 	&	&\rightarrow& 1 & \leftarrow&	&	\\
1	&\cdots	& 0	& \uparrow & 0	&\cdots & 0	\\
	&	&\vdots	&  	&\vdots &	&	\\
	&	& 0	&	& 1	&	&
\end{array}
\end{displaymath}
But the next figure shows that the tilings that result from
looking at the column and the row don't match, so we have a contradiction.
\begin{center}\includegraphics{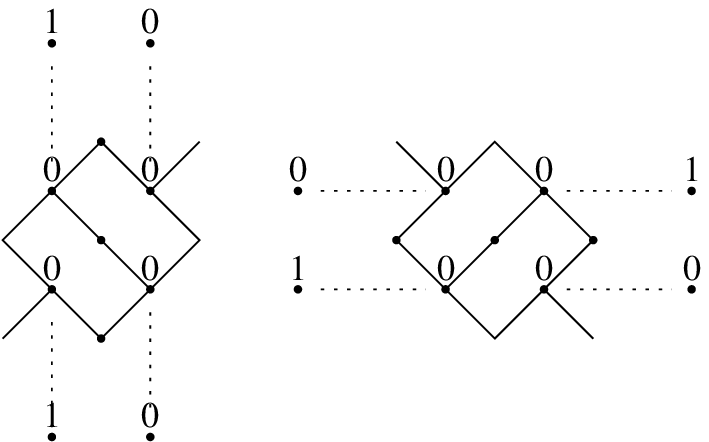}\end{center}

Let us now assume that $\Box =0$, and find a contradiction. If we have a 0
in $A$, then the arrows on either side point in the same
direction.  This is true for the arrows in the row and in the column.  In
the next figure, we have chosen to make the arrows point up and to the
left.  Because of symmetry, we will get the same result for any
equivalent configuration. 
\begin{displaymath}
\begin{array}{ccccccc}
	&	& 1	&	& 0	&	&	\\
	&	&\vdots	&	&\vdots	&	&	\\
0	&\cdots	& 0	& \uparrow & 0	&\cdots & 1	\\
 	&	&\leftarrow & 0 & \leftarrow&	&	\\
1	&\cdots	& 0	& \uparrow& 0	&\cdots & 0	\\
	&	&\vdots	&  	&\vdots &	&	\\
	&	& 0	&	& 1	&	&
\end{array}
\end{displaymath}
But the next figure shows that the tilings that result from
looking at the column and the row don't match, so we have a 
contradiction.
\begin{center}\includegraphics{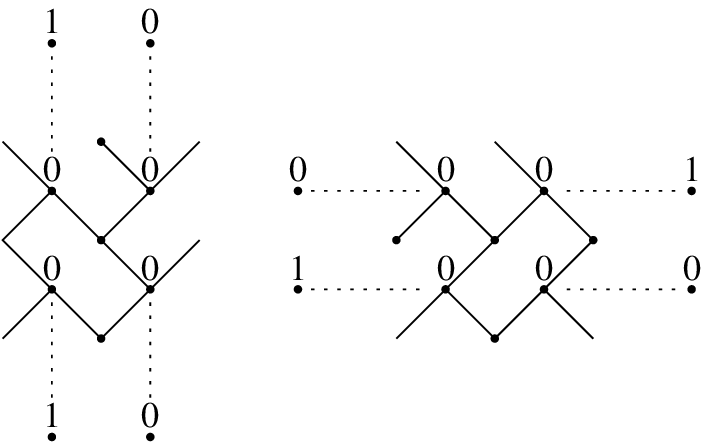}\end{center}
This concludes our proof of Lemma \ref{alg}.

\section{Baxter Permutations and $-1$'s}\label{proof}

With an algorithm in hand that will produce the unique order $n$ ASM
that is compatible with an order $n+1$ permutation matrix, we are
prepared to show that the conditions in the permutation that produce
$-1$'s in the smaller ASM are identical to the conditions that would
prevent that permutation from being Baxter.

First, assume that a permutation $B$ is compatible with a smaller ASM
$A$ that contains no $-1$'s.  If $A$ it contains no $-1$'s then
each  row of $A$ must contain exactly one $\Box$.  Remember that
the $\Box$ represents the place where the direction that the 1's of
$B$ are located switches from up to down or down to up.  Then for
each row, the location of the $\Box$ is the natural place to draw a
dividing line between the columns to pass the Baxter test.
Consequently,  $B$ is Baxter.

Assume that a permutation $B$ is compatible with a smaller ASM
$A$ that contains at least one $-1$.  The row of $A$ that contains the
$-1$ has at least three $\Box$s, and consequently, has 
no natural place to draw the dividing line between columns.  Let us
show that this will keep $B$ from being Baxter.
From our earlier discussion, it is apparent that $B$ contains a
windmilled configuration of 1's. 
\begin{displaymath}
B = 
\left(
\begin{array}{cc cc cc cc}
      	&	&	&\vdots	&\vdots	&	&	&	\\
      	&	&	& 1	& 0	&	&	& 	\\
      	&	&	&\vdots	&\vdots	&	&	&	\\
\cdots	& 0	&\cdots	& 0	& 0	&\cdots & 1	&\cdots	\\
\cdots	& 1	&\cdots	& 0	& 0	&\cdots & 0	&\cdots	\\
      	&	&	&\vdots	&\vdots	&	&	&	\\
      	&	&	& 0	& 1	&	&	&	\\
      	&	&	&\vdots	&\vdots	&	&	&	\\
\end{array}
\right)
\end{displaymath}
It is easily checked that a permutation that looks like this is not
Baxter, because there can be no vertical dividing line between columns
that properly segregates the 1's into two groups.

This concludes the proof of Theorem \ref{bpth}.

\section{Thanks}

I would like to thank James Propp for counting pairs of compatible
permutation matrices for small values of $n$, for making the
conjecture about the connection with Baxter permutations, and for
bringing the problem to my attention.  The other members of the
Spatial Systems Lab at the University of Wisconsin-Madison also
deserve thanks for their help and support with this proof and with
this paper, especially Dominic Johann-Berkel for finding out what a
Baxter permutation was.

This work was supported by the NSF through its VIGRE program, as
administered by the University of Wisconsin.


\begin{thebibliography}{9}

 \bibitem{baxter}
{\sc G. Baxter.}
 ``On fixed points of the composite of commuting functions.''
\emph{Proceedings of the American Mathematical Society} 15 (1964) 851-855.

\bibitem{dulucq}
{\sc S.~Dulucq and O.~Guibert.}
``Baxter Permutations.''
\emph{Discrete Mathematics} 180 (1998) 143-156.

\bibitem{EKLP}
{\sc N.~Elkies, G.~Kuperberg, M.~Larsen, and J.~Propp}.
``Alter\-nating-Sign Matrices and Domino Tilings (Part I).''
\emph{Journal of Algebraic Combinatorics} 1 (1992) 111-132.

\bibitem{kuo}
{\sc Eric H. Kuo.}
``Applications of Graphical Condensation for Enumerating Matchings and
    Tilings.'' 
\texttt{http://arxiv.org/abs/math.CO/0304090} (2003).

\bibitem{robbins}
{\sc D.~P.~Robbins and H.~Rumsey, Jr.}
``Determinants and Alternating Sign Matrices.''
\emph{Advances in Mathematics} 62 (1986) 169-184.

\bibitem{count}
{\sc F.~R.~K.~Chung, R.~L.~Graham, V.~E.~Hoggatt, Jr., and M.~Kleiman.}
``The Number of Baxter Permutations.''
\emph{Journal of Combinatorial Theory}, Series A 24 (1978) 382-394.

\end{thebibliography}
\end{document}